# A Versatile Stochastic Duel Game

SONG-KYOO (AMANG) KIM

**ABSTRACT**

This paper deals with a standard stochastic game model with a continuum of states under the duel-type setup. It newly proposes a hybrid model of game theory and the fluctuation process, which could be applied for various practical decision making situations. The unique theoretical stochastic game model is targeted to analyze a two-person duel-type game in the time domain. The parameters for strategic decisions including the moments of crossings, prior crossings, and the optimal number of iterations to get the highest winning chance are obtained by the compact closed joint functional. This paper also demonstrates the usage of a new time based stochastic game model by analyzing a conventional duel game model in the distance domain and briefly explains how to build strategies for an atypical business case to show how this theoretical model works.

## I. INTRODUCTION

Game theory has been applied for various strategic situations and also developed to solve real-world issues innovatively [1-3]. Specifically, a conventional duel game is an arranged engagement in a combat situation between two players, with matched weapons in accordance with the agreed rules under different conditions [4]. A conventional duel game model accurately describes the conditions in the distance domain and finds when and who could win the battle even at the beginning. Unlike conventional duel games, the duel game in this paper deals with shooting a single target rather than shooting each other. Hence, each player could choose either to "shoot" or "wait" for one step closer to the target during his/her turn (or iteration). The backward induction provides a simple solution: regardless of whether you have better or a worse shot, the shooting moment when the sum of the success probabilities passes the threshold is the most critical [4]. On the other hand, this proposed game is a variant of a real-time based antagonistic stochastic game, and it adapts a duel game on top of conventional antagonistic stochastic games [5-9]. Each player at random times with random impacts can take the best shooting after passing the fixed threshold. According to the backward induction, a player has the chance to take a shot when his/her underlying threshold is crossed. Upon that time (referred to as the first passage time),the player can perform successful shooting [10]. This kind of research has been massively studied in the marketing area [11-17], but none of the research has studied mathematical approaches. Although the new mathematical game model in this paper is relatively restricted, this newly proposed model could be applied to various business decision making issues, especially strategic marketing decisions for new product launch.

In the stochastic duel model, the actions of players are formalized by marked point processes to identify the chances for the successful shootings at certain points in the time domain. The processes evolve until either one of the processes crosses its fixed threshold of success probabilities. Once a threshold is reached at some point in time, the associated player has the highest chance to win the game. This standard stopping game model could be applied to various business decision making situations, although the rules and conditions are relatively restricted. Various duel-type games have been studied since the 1970s [4, 18-21], and the key lesson from them is that the best decision that matters is when to do, rather than what to do [4]. Unfortunately, conventional duel-type games only consider the deterministic turn around (iteration) for shooting by the success probabilities based on the distance between players. However, the proposed model in the paper considers random iterations for decision making. The strategic decision making situation in the smartphone business could adopt a versatile antagonistic game. Although there are many smartphone manufacturers in the world, Apple and Samsung are the most dominant manufacturers on the market. This game model is designed as two players (i.e., Apple and Samsung). Because smartphone technologies are rapidly changing, devices need to be upgraded even after launching products. Strategically, a phone manufacturer that implements more technologies (or features) has a greater chance of winning, but it requires more time and resources for research and development activities. More importantly, a firm could lose its market share if the company launches a product that fails to satisfy customers. Therefore, an initial product should be appealing enough to dominate the market. However, at the same time, the firm could also fail if a product is launched too late compared to rivals.

The article presents a versatile stochastic duel game with complete information, which means both players know the success probabilities in the time domain. Unlike conventional (in the distance domain) duel games, each player might not have the same iteration periods, and the iteration periods for each round might be different even within the same player (i.e.,stochastic). Lastly, the paper demonstrates the special case of the stochastic duel game, which is based on the deterministic same iteration times for both players. This special case shows how the time domain duel game model is related to the distance based duel games.

## II. ANTAGONISTIC STOCHASTIC DUEL GAME

The antagonistic duel game of two players (called "A" and "B") are introduced and both players know the completed information regarding the success probabilities based on the time. Each player has two strategies either "shoot" or "wait" and choose one strategy at the certain points of time. Let $A(s)$ be player A's payoff (related) function based on the continues time $s$ and $B(t)$ be player B's payoff (related) function based on the continues time $t$. Both functions are assigned as follows:

$$\{A(s) : 0 \leq A(s) \leq A(s + \Delta),\ s \in [0, s_{\max}), s_{\max} \in \mathbb{R}_+,\ \Delta > 0\}, \tag{1}$$

$$\{B(t) : 0 \leq B(t) \leq B(t + \Delta),\ t \in [0, t_{\max}), t_{\max} \in \mathbb{R}_+,\ \Delta > 0\}. \tag{2}$$

The accumulative success probability functions of both players could be as follows:

$$P_a(s) := \frac{A(s)}{A(s_{\max})},\ P_b(t) := \frac{B(s)}{B(s_{\max})}. \tag{3}$$

The strategic decision for a duel game is to find the moment when a player will have the best chance to hit the other. Both accumulative success probabilities could be any arbitrary incremental continuous function which reaches 1 when the time $s$ (or $t$) goes to the allowed maximum ($s_{\max}$) and, typically, it goes to the infinity. In the duel game, there is a certain point that maximizes the chance for succeeding the shoot (i.e., success probability). This optimal point becomes the moment of the success in the continuous time domain and this moment $t^*$ is defined as follows:

$$t^* = \inf\{t \geq 0 : P_a(t) + P_b(t) \geq 1\}. \tag{4}$$

It is noted that each player can make the decision at the certain points of the time and that is the reason why it becomes a discrete time series even though the success probabilities are continuous functions. Let $(\Omega, \mathcal{F}(\Omega), P)$ be a probability time space and let $\mathcal{F}_S, \mathcal{F}_T \subseteq \mathcal{F}(\Omega)$ be independent $\sigma$-subalgebras. Suppose

$$S := \sum_{j \geq 0} \varepsilon_{S_j},\ 0 = S_0 < S_1 < \ldots,\ \text{a.s.} \tag{5}$$

$$T := \sum_{k \geq 0} \varepsilon_{T_j},\ 0 = T_0 < T_1 < \ldots,\ \text{a.s.} \tag{6}$$

are $\mathcal{F}_S$-measurable and $\mathcal{F}_T$-measurable renewal point processes with the following notation:

$$\sigma := \begin{cases} \sigma_j = S_j - S_{j-1}, & j = 1, 2, \ldots, \\ 0, & j \leq 0, \end{cases} \tag{7}$$

$$\tau := \begin{cases} \tau_k = T_k - T_{k-1}, & k = 1, 2, \ldots, \\ 0, & k \leq 0. \end{cases} \tag{8}$$

The game in this case is a stochastic process describing the evolution of a conflict between players A and player B based on the fully known information (i.e., the success probabilities of players). Only on the $j$-th epoch $S_j$, player A could make the decision either for taking a shoot or for waiting until another turn (iteration) $S_{j+1}$. He will have the best chance to hit the player B exceeds its respective threshold $U$ (or $V$ for player B). To further formalize the game, the exit indices are introduced as follows:

$$\mu := \inf\{j : S_j = \sigma_0(\leq 0) + \sigma_1 + \cdots + \sigma_j \geq U\} \tag{9}$$

and

$$\nu := \inf\{k : T_k = \tau_0(= 0) + \tau_1 + \cdots + \tau_k \geq V\} \tag{10}$$

and $\sigma_0 < 0$ indicates that player A is starting first in this game. In the case of the duel type games, the threshold of each player could be converged into one value $t^*$ which will be introduced later. Player A will have the best chance to succeed for shooting compare to the failure chance of player B ($P_a(S_\mu)$ and $1 - P_b(T_\nu)$ respectively). Hence, player A has the highest success probability of shooting at time $S_\mu$, unless player B does not reach his best shooting chance at time $T_\nu$. Thus, the game is ended at $min\{S_\mu, T_\nu\}$. However, we are targeting the confined duel game for player A on trace $\sigma$-algebra $\mathcal{F}(\Omega) \cap \{P_a(S_\mu) + P_b(T_\nu) \geq 1\} \cap \{S_\mu \leq T_\nu\}$ (i. e., the game with player A obtains the best chance for shooting first). The first passage time $S_\mu$ is the associated time from the confined game. The functional

$$\begin{aligned} \Phi_{\mu\nu} &= \Phi_{\mu\nu}(\theta_0, \theta_1, \vartheta_0, \vartheta_1) \\ &= \mathbb{E}\left[e^{-\theta_0 S_{\mu-1} - \theta_1 S_\mu} \cdot e^{-\vartheta_0 T_{\nu-1} - \vartheta_1 T_\nu} \cdot \mathbf{1}_{\{S_\mu \leq T_\nu\}} \mathbf{1}_{\{A(S_\mu) + B(T_\nu) \geq 1\}}\right] \end{aligned} \quad (11)$$

of the game will represent the status of both players upon *exit time* $S_\mu$ and *pre-exit time* $S_{\mu-1}$ (Dshalalow, 1995). The latter is of particular interest, because player A wants to predict not only his time of the highest chance, but also the moment of the next highest chance prior to this. The Theorem 1 establishes an explicit formula for $\Phi_{\mu\nu}$ and we abbreviate with (12)-(19):

$$\gamma(x, t) = e^{-xt}, \quad (12)$$

$$\gamma_0(t) := \gamma(v, t), \quad (13)$$
$$\gamma_1(t) := \gamma(\vartheta_0 + v, t), \quad (14)$$
$$\gamma_2(t) := \gamma(\vartheta_0 + \vartheta_1 + v, t), \quad (15)$$

$$\Gamma_0(t) := \gamma(u, t), \quad (16)$$
$$\Gamma_1(t) := \gamma(\theta_0 + u, t), \quad (17)$$
$$\Gamma_2(t) := \gamma(\theta_0 + \theta_1 + u, t), \quad (18)$$
$$\Gamma(t) := \gamma_2(t) \cdot \Gamma_2(t). \quad (19)$$

We will use the Laplace-Carson transform

$$\widehat{\mathcal{L}}_{pq}(\bullet)(u, v) = uv \int_{p=0}^{\infty} \int_{p=0}^{\infty} e^{-up-vq}(\bullet) d(p, q), \; \operatorname{Re}(u) > 0, \; \operatorname{Re}(v) > 0, \quad (20)$$

with the inverse

$$\widehat{\mathcal{L}}_{uv}^{-1}(\bullet)(p, q) = \mathcal{L}^{-1}\left(\bullet \frac{1}{uv}\right), \quad (21)$$

where $\mathcal{L}^{-1}$ is the inverse of the bivariate Laplace transform [23].

**Theorem 1.** The functional $\Phi_{\mu\nu}$ of the game on trace race $\sigma$-algebra $\mathcal{F}(\Omega) \cap \{A(s_\mu) + B(s_\mu) \geq 1\} \cap \{S_\mu \leq T_\nu\}$ satisfies the following formula:

$$\Phi_{\mu\nu} = \widehat{\mathcal{L}}_{uv}^{-1}\left(\mathbb{E}\left[\frac{(1 - \gamma_0(\tau))(1 - \Gamma_0(\sigma))}{\gamma_1(\tau)(1 - \gamma_2(\tau))\Gamma_1(\sigma)(1 - \Gamma(\sigma))}\right]\Gamma(t^*)\right)(t^*, t^*). \quad (22)$$

∎

The functional $\Phi_{\mu\nu}$ contains all decision making parameters regarding this game. The information includes the best moments of shooting ($S_\mu, T_\nu$; *exit time*) and the one step before the best moments ($S_{\mu-1}, T_{\nu-1}$; *pre-exit time*) and the optimal number of iterations for both players. The information for both players from the closed functional are as follows:

(for player A)
$$\mathbb{E}[S_\mu] = \lim_{\theta \to 0}\left(-\frac{\partial}{\partial \theta}\right)\Phi_{\mu\nu}(0, \theta, 0, 0), \quad (23)$$

$$\mathbb{E}[S_{\mu-1}] = \lim_{\theta \to 0} \left( -\tfrac{\partial}{\partial \theta} \right) \Phi_{\mu\nu}(\theta, 0, 0, 0), \tag{24}$$

(for player B)
$$\mathbb{E}[T_\nu] = \lim_{\vartheta \to 0} \left( -\tfrac{\partial}{\partial \theta} \right) \Phi_{\mu\nu}(0, 0, 0, \vartheta), \tag{25}$$

$$\mathbb{E}[T_{\nu-1}] = \lim_{\vartheta \to 0} \left( -\tfrac{\partial}{\partial \theta} \right) \Phi_{\mu\nu}(0, 0, \vartheta, 0), \tag{26}$$

and

$$\mu = \left\lfloor \frac{\mathbb{E}[S_\mu]}{\mathbb{E}[\sigma]} \right\rfloor, \quad \nu = \left\lfloor \frac{\mathbb{E}[T_\nu]}{\mathbb{E}[\tau]} \right\rfloor. \tag{27}$$

### III. CASE PRACTICE: NEW PRODUCT LAUNCHING STRATEGY

This section demonstrates a practical application for developing marketing strategy and this case practice shows how to applied the proposed versatile stochastic duel game. The numbers and names are fictional but realistic for this demonstration. Two major smartphone companies are considered and let us make some stories as follows: Samsung just finished developing the flagship smartphone for this year. Samsung should decide to either "Launch" the smartphone, or "Wait" for developing new additional features. One development cycle takes approximately 6 months. Once a development cycle is started, they cannot launch the product until it is completed. In other words, Samsung must wait 6 months for his next opportunity to make another decision once they decide to "Wait." Samsung knows that Apple completes the flagship smartphone 5 months later, and it takes 4 months to add new additional features by Apple after completing the initial product. Since both companies have ample experience in the smartphone industry, they know the probability of success for either choice. Samsung wants to know the best time to launch his new flagship smartphone. If the flagship product is launched too early, consumers would not choose the product because it does not have enough new technologies to be attractive. However, Samsung may also lose the game if Apple launches their flagship device before Samsung launches. This versatile duel game could be easily demonstrated to provide the best strategies for Samsung (and Apple). It is noted that this section is designed only for showing how the mathematical model could be applied in the marketing decision problem. Hence, this section only provides the condensed results. The development cycles of both companies are assumed to be exponentially distributed and the parameters based on the current scenario are as follows:

| Parameter | Value | Description |
|---|---|---|
| $\widetilde{\delta_0}\ (= \mathbb{E}[S_0])$ | 0 | Initial development cycle of Samsung |
| $\widetilde{\delta}\ (= \mathbb{E}[\sigma])$ | 6 | Development cycle of Samsung |
| $\mathbb{E}[T_0]$ | 5 | Initial development cycle of Apple (after Samsung) |
| $\mathbb{E}[\tau]$ | 4 | Development cycle of Apple |

[**Table 1.** Initial setup values for the case practice]

It is noted that some parameters of this setup, including the development cycles, are not carefully evaluated, but are realistic because these are chosen based on information from experts in the smartphone industry. Once all the parameters are mapped, we can get the decision parameters from the simulated results.

| Parameters | Values | Descriptions |
|---|---|---|
| $t^*$ | 17.95 | Best moment that whoever player win the game |
| $\mathbb{E}[\mu]$ | 3 | Number of development cycle that is best for Samsung |
| $\mathbb{E}[S_\mu]$ | 18 | Best time to launch the product for Samsung |
| $\mathbb{E}[\nu]$ | 4 | Number of development cycle that is best for Apple |
| $\mathbb{E}[T_\nu]$ | 21 | Best time to launch the product for Apple |

[**Table 2.** Decision making results]

According to the simulated results, Samsung has the opportunity to launch the flagship product after passing the moment $t^*$ but before launching the product by Apple. But Samsung should run at least 3 development cycles to add new features before launching (i.e., $\mu = 3$) after first commercially ready product is developed ( $= S_0$). If both companies launch their flagship smartphones before 18 months ( $< t^*$), customers would not buy their products because they lack attractive features. In the setup of above scenario, Samsung is mostly likely win the game because the best moment of shooting by Samsung ( $= S_\mu$) is closer to the threshold $t^*$ than what Apple has (i. e., $t^* < S_\mu < T_\nu$). In this particular scenario, releasing the smartphone on $T_\nu$ by Apple is no longer the optimal strategy. Therefore, it might be better that Apple takes the risk to launch the product at $T_{\nu-1}$ ( $= 17$ months) instead of $T_\nu$ ( $= 21$ months) because Apple will lose this game when Samsung launches the flagship product at $S_\mu$. Apple will probabilistically lose this game but it is noted that the solution of the game is not deterministic. In other words, a player with a higher probability of success is not guaranteed to win a game.

## IV. Conclusion

A new type of antagonistic stochastic duel game was studied. In this versatile stochastic duel game, both players could have random iterations in the time domain. A joint functional of the standard stopping game was constructed to analyze the decision making parameters, which indicated the best winning strategies of players in the time domain stochastic game. Compact closed forms from the Laplace-Carson transforms were obtained, and the hybrid stochastic game model was more flexible to solve various duel-type game problems more effectively. The analytical approach by applying the proposed model to a conventional duel game was fully described for better understanding the core of the model. Furthermore, an actual application in the smartphone market was provided for the direct implementation of a versatile stochastic duel game in real-world decision making situations. An extended version of the duel-type game for multiple players is currently under development as a successor of this research.